\documentclass[12pt,leqno]{article}
\usepackage{amstext}
\usepackage{theorem}
\usepackage{amssymb}
\usepackage{latexsym}
\newtheorem{thm}{Theorem.}[section]
\newtheorem{prop}[thm]{Proposition.}

\newtheorem{lem}[thm]{Lemma.}
\newtheorem{rem}[thm]{Remark.}

\newtheorem{prop-def}[thm]{Proposition-Definition.}

\def\pmb#1{\setbox0=\hbox{#1}%
\kern-.025em\copy0\kern-\wd0
\kern.05em\copy0\kern-\wd0
\kern-.025em\raise.0433em\box0}

\newcommand{\be}{\begin{enumerate}}
\newcommand{\ee}{\end{enumerate}}

\newcommand{\br}{\begin{array}}
\newcommand{\er}{\end{array}}

\newfont{\ninemsbm}{msbm10 scaled 0900}
\newfont{\tenmsbm}{msbm10 scaled 1100}
\newfont{\nineeufb}{eufb10 scaled 0900}
\newfont{\teneufb}{eufb10 scaled 1100}
\newfont{\teneusm}{eusm10 scaled 1100}
\newfont{\nineeusm}{eusm10 scaled 0900}
\newfont{\bal}{cmmib10 scaled 0900}
\newfont{\tencmmib}{cmmib10 scaled 1000}

\newcommand{\gotm}[1]{\mbox{{\teneufb {#1}}}}

\newcommand{\PP}{\mbox{$I\!\!P$}}

   \font\tenmsb=msbm10 scaled\magstep 1        
   \font\sevenmsb=msbm7 scaled \magstep 1
   \font\fivemsb=msbm5 scaled \magstep 1
\newfam\msbfam
      \textfont\msbfam=\tenmsb 
      \scriptfont\msbfam=\sevenmsb
      \scriptscriptfont\msbfam=\fivemsb

\font\tenblah=msbm10 scaled\magstep 3        
   \font\sevenblah=msbm7 scaled \magstep 1
   \font\fiveblah=msbm5 scaled \magstep 1
\newfam\blahfam
      \textfont\blahfam=\tenblah 
      \scriptfont\blahfam=\sevenblah
      \scriptscriptfont\blahfam=\fiveblah

\begin{document}
\begin{center}{\huge\bf New evidence for Green's conjecture on
   syzygies of canonical curves
\vspace{3ex}}\\
A. Hirschowitz\hspace{3cm}S. Ramanan\end{center}
\section{Introduction}
Some twelve years ago, Mark Green [G] made a few conjectures
regarding the behaviour of syzygies of a curve $C$ imbedded in
$\PP^n$ by a complete linear system. The so-called {\it generic} Green
conjecture on canonical curves pertains to this question when
the linear system is the canonical one and the curve is generic
in the moduli, and predicts what are the numbers of syzygies in that
case. Green and Lazarsfeld [GL] have observed that curves with 
nonmaximal Clifford index have extra syzygies and 
we will call 
{\it specific} Green
conjecture on canonical curves the stronger prediction that the
curves which have the numbers of 
syzygies expected for generic curves 
are precisely those with maximal Clifford index  ($[(g-1)/2]$).
(As a matter of fact, the full  Green
conjecture on canonical curves relates more closely the Clifford
number
with the existence of extra syzygies.)
Many attempts have been made to settle
this question, and some nice results have been obtained ([Sch][V]).

In this note, we work over an algebraically closed field of
arbitrary
characteristic and prove that, as stated above, the generic and specific
Green conjectures for canonical curves are equivalent at least when
the genus $g$ is odd. 
Let $C$ be a  curve canonically imbedded in $\PP^{g-1}$
with ideal sheaf
 ${\cal
   I}_C$. We denote by $Q$ the universal quotient on
$\PP^{g-1}$, so that $Q(1)$ is the tangent bundle, and by $Q_C$ its
restriction to $C$. 
It is generally known (see [P-R] for example)
that extra syzygies appear when,
for some $i\leq [(g-1)/2]$, the natural map
$\Lambda ^i (\Gamma (C,Q_C)) \rightarrow \Gamma (C, \Lambda
^iQ_C)$ is not surjective. It is easy to see that the relevant quotient
of $\Gamma (C, \Lambda
^iQ_C)$ by $\Lambda ^i (\Gamma (C,Q_C))$ is isomorphic to 
 $\Gamma (\Lambda^{i+1}Q \otimes {\cal I}_C(1))$ (cf 2.1).

\begin{thm} Let $g = 2k-1 \geq 5$ be an odd integer. If the generic curve $C$
  of
  genus $g$ has the expected number of syzygies (i.e. 
$\Gamma (\Lambda ^{k} Q \otimes {\cal I}_C(1))=0$),
then so does any curve of genus $g$ with maximal Clifford index, namely 
$k-1.$
\end{thm}

To prove this,
we compute
a virtual (divisor) class $v$ for the
locus (in the moduli) of curves $C$ for which
the cohomology group
$\Gamma (\Lambda^{k} Q \otimes {\cal I}_C(1))$ does not vanish.
Once $v$ is computed, we compare it with
the class $c$ of
the locus of $k$-gonal curves (these are 
curves with non-maximal Clifford index in our
case), which, 
thanks to Harris and Mumford [HM], is already known,
and we find that $v=(k-1)c$. We conclude by proving that the 
generic $k$-gonal curve has at least $k-1$ extra syzygies, which 
implies that the $k$-gonal locus 
occurs with multiplicity
$k-1$ in $v$,
leaving
no room for another component. Our proof gives another
consequence
of the generic Green's conjecture, namely that the number of extra
syzygies (more precisely  
$h^0(\Lambda^{k} Q \otimes {\cal I}_C(1))$)
is exactly $k-1$ for any $k$-gonal curve $C$ in the smooth part of the
$k$-gonal locus.
Finally, we observe
that our argument fails completely in the case of even genus, where the
expected codimension of our jump locus is no more one.
       
\section { Preliminaries on syzygies} \label {prel}

We collect here a few useful remarks on syzygies.

\begin{prop} Let $S$ be a linearly normal
subscheme of $\PP^n $
(i.e.  $\Gamma (\PP^n,{\cal O}(1)) \to
\Gamma (S, {\cal O}(1))$ is an isomorphism) with ideal sheaf ${\cal
  I}_S$.
Then the cokernel of
$$
\Lambda^ i\Gamma (\PP^n,Q) = \Gamma (\PP^n, \Lambda ^i Q)\rightarrow
\Gamma (S,\Lambda^i Q_S) 
$$ 
is canonically isomorphic to 
$$\Gamma (\PP^n,\Lambda^{i+1} Q\otimes {\cal I}_S(1)).$$
\end{prop}

\noindent{\bf Proof.}

Consider the following exact sequence of sheaves on $P(V) =
\PP^n$:

$$
0\rightarrow {\cal O}(-1) \rightarrow V_P \rightarrow Q
\rightarrow 0
$$

By taking the exterior $(i+1)$-th power and tensoring with ${\cal
O}(1)$, we get the exact sequence:

$$0\rightarrow \Lambda^{i}Q \rightarrow \Lambda^{i+1} V_P (1)\rightarrow
\Lambda^ {i+1}Q (1) \rightarrow 0.
$$

This exact sequence of vector bundles remains
exact on tensorisation by ${\cal I}_S$ as well as ${\cal O}_S$.
Thus we get the commutative diagram

$$\begin{array}{ccccccccc}
  &    &    0        & &     0        & &     0     & &\\
  
  &    & \downarrow  & & \downarrow  & & \downarrow & &\\
  
0 & \rightarrow & \Lambda^i Q \otimes {\cal I}_S & \rightarrow &
\Lambda^{i+1} V_P \otimes {\cal I}_S(1) & \rightarrow & \Lambda^{i+1} Q
\otimes {\cal I}_S(1)  & \rightarrow & 0 \\

  &    & \downarrow  & & \downarrow  & & \downarrow & &\\

0 & \rightarrow & \Lambda^i Q & \rightarrow & \Lambda^{i+1} V_P (1) &
\rightarrow & \Lambda^{i+1} Q(1) & \rightarrow & 0\\

  &    & \downarrow  & & \downarrow  & & \downarrow & &\\

0 & \rightarrow & \Lambda^i Q_S & \rightarrow & \Lambda^{i+1} V_S(1)
 & \rightarrow & \Lambda^{i+1} Q_S(1) & \rightarrow & 0 \\ 

  &    & \downarrow  & & \downarrow  & & \downarrow & &\\

  &    &    0        & &      0      & &     0      & &  

\end{array}
$$

Now apply the section functor $\Gamma $ : the middle row remains
exact. Thus we may apply the snake lemma to the two lower rows. This
yields the desired isomorphism because under our
assumption, $\Gamma (\Lambda ^{i+1}V_P(1)) \to \Gamma (\Lambda
^{i+1}V_S(1))$ is an isomorphism as well. $\hfill\square$

\begin{rem}
Thus we will think of $\Gamma (\Lambda ^j Q \otimes {\cal I}_C(1))$
as the space of extra syzygies. From this point of view,
extra syzygies behave in a monotonic way with respect to the
degree $j$ and the subvariety $C$:

a) If $C \subset S$ are two subvarieties of $\PP^{g-1}$,
then  $h^0 (\Lambda ^j Q \otimes {\cal I}_C(1))\geq 
h^0 (\Lambda ^j Q \otimes {\cal I}_S(1))$. We will estimate
syzygies of our canonical curves by using a scroll $S$ containing
them.

b) If $i < j$, then  $h^0 (\Lambda ^j Q \otimes {\cal I}_C(1))\geq 
h^0 (\Lambda ^i Q \otimes {\cal I}_C(1)).$

The above proposition is applicable with our canonical curve: $S = C$.
Also, the Clifford index of the  generic curve of genus $g$, is well-known
to be $[(g-1)/2]$. Finally, if $ i < j$ then
$\Gamma (\Lambda ^j Q \otimes {\cal I}_C(1))=0$ implies
$\Gamma (\Lambda ^i Q \otimes {\cal I}_C(1))=0$
so that we have an equivalent formulation of the
specific conjecture of Green:

\noindent{\bf Specific Green's conjecture.} {\it Let $C$ be a
canonically imbedded curve with maximum Clifford index 
$[(g - 1)/2]$. Then
$\Gamma (\Lambda ^j Q \otimes {\cal I}_C(1))$ is zero for
$j = [(g+1)/2]$.}
\end {rem}

We will use in  Section \ref {scroll} the following semi-continuity statement:

\begin{prop} \label {scs} Let $p: W \rightarrow T$ 
be a smooth family of projective
  varieties parametrized by  the spectrum $T$ of a discrete valuation
  ring. We suppose that $W$ is 
endowed with a line bundle ${\cal L}$,
 that $h^0(W_t, {\cal L}_t)$ is constant and that
for each point $t \in T$, ${\cal L}_t$ is generated by global
  sections. This
yields a $T$-morphism $m$ from $W$ to $\PP(p_*{\cal L})$. We denote 
by $I_t$
the ideal sheaf of $m(W_t)$ and by $Q_t$ the tautological quotient
  bundle 
on the fibre $\PP_t=\PP(H^0(W_t,  {\cal
  L}_t))$.
Then for any $i$, the dimension 
$h^0 (\PP_t, \Lambda^iQ_t \otimes I_t(1))$ is upper-semi-continuous.
\end{prop}

\noindent {Proof.}
By properness of the Hilbert scheme, there exists
a $T$-flat subscheme $\bar W$ of $\PP(p_*{\cal L})$ with the property
that
its fibre over the general point $t_1$ of $T$ is $m(W_{t_1})$. By
continuity,
its special fibre $\bar W_{t_0}$ contains $m(W_{t_0})$ (indeed, they
are equal, but we don't need this).
Thus, by inclusion, we have
$$h^0 (\PP_{t_0}, \Lambda^iQ_{t_0} \otimes I_{t_0}(1)) \geq
h^0 (\PP_{t_0}, \Lambda^iQ_{t_0} \otimes I_{\bar W_{t_0}}(1)),$$
and by semi-continuity,
$$h^0 (\PP_{t_0}, \Lambda^iQ_{t_0} \otimes I_{\bar W_{t_0}}(1)) \geq
h^0 (\PP_{t_1}, \Lambda^iQ_{t_1} \otimes I_{\bar W_{t_1}}(1)),$$
which altogether prove our claim.$\hfill \square$

\section{The syzygy locus in the case of odd genus} \label {S3}

In this section, we write ${\cal M}={\cal M}^{o}_g$ for 
the open subvariety of ${\cal M}_g$ consisting
of points that represent isomorphism classes of smooth curves with
trivial automorphism group. What we need to know of ${\cal M}$ is
that an effective divisor on it which is rationally equivalent to zero
is indeed zero: this is for instance because ${\cal M}$ has a
projective compactification with two-codimensional boundary (cf  e.g. [A]).

Let $x$ be  a point in ${\cal M}$ and $C$ the corresponding curve.
We consider the canonical imbedding of $C$ in $\PP^{g-1}$
($C$ is not hyperelliptic);
we denote by ${\cal I}_C$ the ideal sheaf of $C$ 
and by $Q$ the tautological quotient bundle of rank $g-1$ on $\PP^{g-1}$. 
Finally, we denote by $S_{g}$ 
the locus in ${\cal M}$ of (points corresponding to) curves
$C$ satisfying $\Gamma (\Lambda^{k} Q \otimes {\cal I}_C(1)) \neq 0.$ 
As a jump locus, $S_g$ has
a natural Cartier divisor structure (see e.g. the proof of the next proposition)
and
we
compare its  class  in the Picard group of ${\cal
  M}$ with
the class $c$ of the $k$-gonal locus (cf  [HM]).

\begin{prop} Let  $g = 2k-1\geq 5$ be an odd
  integer
such that the generic curve $C$ of genus $g$ satisfies
$\Gamma (\Lambda^{k} Q \otimes {\cal I}_C(1))=0.$ Then, in the Picard group of
${\cal M}$,
 the rational class $v$ of $S_{g}$ is $(k-1)c$.
\end{prop}
\noindent{\bf Proof.}

 There exists a universal curve ${\cal C}$ over ${\cal M}$, 
that is to say a smooth variety $\cal C$ and a smooth
projective morphism $\pi :\cal C\to \cal M$ such that for any
$x\in \cal M$ the fibre of $\pi $ over $x$ is the curve of
genus $g$ whose isomorphism class is given by the point $x$. Let
$\omega = \omega_{\pi }$ be the cotangent bundle along the
fibres and $E$ its direct image on ${\cal M}$ by $\pi $. Then $\pi$
factors
through
the natural canonical imbedding of ${\cal C}$ in the projective
bundle $p:\PP = \PP (E) \to {\cal M}$. Let ${\cal I}$ be the ideal
sheaf
of ${\cal C}$ in $\PP$.
The relatively ample (hyperplane) line bundle along the fibres
of $\PP $ will be denoted as usual by $ {\cal O}_p(1)$.
Finally we write again $Q$ for the vector bundle on $\PP $ given by the
exact sequence

$$
0\rightarrow {\cal O}_p(-1) \rightarrow p ^*(E)^*
\rightarrow Q \rightarrow 0.
$$
 
Observe that $p _*(\Lambda ^l Q(1))$ is a vector bundle of rank
${g\choose l}g - {g\choose l - 1}$. 
Similarly, on each fibre ${\cal C}_x$, $\Lambda ^l Q_{{\cal C}_x}\otimes
\omega_{{\cal C}_x}$
is semi-stable (cf [PR]) of slope $2l + 2g -
2$, thus non-special, and 
$p _*(\Lambda ^l Q(1)\otimes {\cal O}_{\cal C})$ is also a vector
bundle, of rank
${g-1\choose l}(2l + g -
1)$, for each $l>0.$
 
Substituting $k$ for $l$, the
 restriction from $\PP$ to the universal curve yields a morphism
$r$ from $p _*(\Lambda ^k Q(1))$ to 
$p _*(\Lambda ^k Q(1)\otimes {\cal O}_{\cal C})$, and our assumption
means that this morphism is injective (at the generic point).
We observe
that the
two vector bundles have the same rank, namely ${2k - 2\choose k
}(4k - 2) = {2k - 1\choose k - 1}(2k - 2)$. Thus the above
map defines a (degeneracy) divisor in ${\cal M}$ and this is $S_g$, 
by definition. Its (virtual) rational class is
$v=c_1(p _*(\Lambda ^k Q(1)\otimes {\cal O}_{\cal C}))-c_1(p _*(\Lambda ^k Q(1))).$

We will
compute this class in $Pic({\cal M})$ as a
multiple of $\lambda= c_1(E)$. We will start with the following
computation in the appropriate Grothendieck group $K$.  Let $t$
be an indeterminate and for any vector bundle $V$, let $\lambda
_t(V)$ denote the element $\sum t^i \Lambda ^i(V)$ in $K[[t]]$.
This extends to a homomorphism of $K$ into the multiplicative
group consisting of power series with constant term 1 in
$K[[t]]$, and this map is still denoted by $\lambda _t$.
Consider now $x = p _!(\lambda _t(Q).{\cal I}(1)) = p
_!(\lambda _t(Q).({\cal O}_p(1) - {\cal O}_{\cal C}(1)))$. Substitute $Q =
p ^*(E^*) - {\cal O}_p(-1)$.  Then we obtain
\begin{eqnarray*}
x & = & p _!( \frac {\lambda _t(p ^*(E^*))}{\lambda _t({\cal
        O}(-1))} ({\cal O}_p(1) - {\cal O}_{\cal C}(1))\cr
  & = & \lambda _t(E^*) p _!
\left ( \frac {{\cal O}_p(1) - {\cal O}_{\cal C}(1)}{1
        + t{\cal O}_p(-1)}\right )\cr
  & = & {\lambda }_t(E^*) \sum _{j = 0} ^{j = \infty } (-1)^j t^j 
        (p _!({\cal O}_p(1-j)) - p _!({\cal O}_{\cal C}(1-j))).
\end{eqnarray*}
Our class $v$ is the coefficient of $t^k$ in the first
Chern class of $-x$.  Observe that $p _! ({\cal O}_p(1 - j)) = 0$,
whenever $2 \leq j \leq g$. Also we have $p _!({\cal O}) = 1$
and $p _!({\cal O}_p(1)) = E. $   On the other hand, $p _!({\cal
O}_{\cal C}(1 - j))$ can be seen to be   $ E^* - 1$ for $j = 0$ and to
be $1 - E$ for $j = 1$. The first Chern class of the direct
images for $j\geq 2$ can be computed by the Grothendieck-Riemann-Roch 
theorem to be $(1 - 6(1 - j) + (1-j)^2) \lambda $ (see
[M]).  The first Chern class of $\lambda _t(E^*)$ is clearly
equal to $-t(1 + t)^{g - 1}\lambda $. Also the rank of $\lambda
_t(E)$ is $(1 + t)^g$, while the rank of $p _!({\cal O}_{\cal C}(n))$
is $(g - 1)(2n - 1)$.  Thus $v$
is equal to $N\lambda $ where $-N$
is the coefficient of $t^k$ in 
$$ 
(1 + t)^g\{ 1 - \sum ^{i =
\infty}_{i = 0} (-1)^i (1 + 6i + 6i^2)t^i \} - t(1 + t)^{g -
1}\{ g - t - (g - 1)\sum _{ i = 0}^{i = \infty } (- 1)^i t^i (1 - 2i)\}. 
$$
  
On the one hand, we have

$(1 + t)^g (1 - \sum _{i = 0} ^{i = \infty } (-1)^i(1 - 6i + 6i^2)t^i)$

$=
(1 + t)^g ( 1 - \sum _{i = 0} ^{i = \infty }(-1)^i(6(i + 1)(i + 2)
- 24(i + 1) + 13)t^i )$

$=(1 + t)^g ( 1 - {12\over (1 + t)^3} + {24\over (1 + t)^2} -
{13\over 1 + t})$

$=(1 + t)^{g - 3}( (1+t)^3 - {13}(1 + t)^2 + {24} (1 + t) - 12)$

$=t(1 + t)^{g - 3} (t^2 - 10t +1),$

and on the other,

$t(1 + t)^{g - 1} ( g - t - (g-1)\sum _{i = 0} ^{i = \infty } (-1)^i(
3 - 2(i + 1))t^i)$\\

$=t(1 + t)^{g - 1}(g - t - (g - 1)({3\over 1 +
t} - {2\over (1 + t)^2}))$\\ 

$=t(1 + t)^{g - 3}((1 + t)^2(g - t) -
(g - 1)3(1 + t) - 2(g - 1))$\\ 

$=t(1 + t)^{g - 3}(-t^3 + (g-2)t^2
+ (-g + 2)t + 1).$

This leads to the determination of $N$ to be the coefficient of
$t^k$ in 
$$t^2(1 +t)^{2k - 4}(-t^2 +(2k - 4)t - (2k - 13)),$$
namely
$$- {2k - 4\choose k - 4} + (2k - 4){2k - 4\choose k - 3} -
(2k - 13){2k - 4\choose k - 2}$$
and this simplifies to 
$$ 6(k +1)(k - 1){(2k - 4)!\over (k-2)!k!}.$$

Now Harris and Mumford [HM] have studied the locus of $k$-gonal
curves in ${\cal M}$ and have shown that this variety
is a divisor whose class is $6(k+1){(2k - 4)!\over (k-2)!k!}\lambda $,
which proves our claim. $\hfill \square$

\section{Syzygies of scrolls} \label {scroll}

Extra syzygies of $k$-gonal curves arise because they lie on scrolls. 
So we start with estimating some syzygies of scrolls.
     
\begin{prop} Let $W$ be a vector bundle on
$\PP ^1$ of rank $k-1$ and degree $k$. We suppose $W$ 
to be globally
generated. We denote by $I_W$ the ideal of the image of the natural morphism 
from $\PP (W)$ into $\PP \Gamma (W)$ and by $Q$ the tautological
quotient
bundle on this projective space. Then
the dimension $h^0 (\PP \Gamma (W), \Lambda ^{k
} Q\otimes I_W(1))$ is at least $k - 1$.
\end{prop}

\noindent {\bf Proof.} By \ref {scs}, we may suppose that $W$ is generic
namely
$W= {\cal O}(1)^{\oplus {k-2}} \oplus {\cal
O}(2)$. In this case,
the natural morphism $\PP(W) \to \PP
\Gamma (W)$ is an imbedding. We will use freely the identification (2.1).
Consider $X = \PP ^1 \times \PP ^1$ and the variety $Y = \PP^1
\times \PP (W)$. Let us denote as usual by $p_1$ and $p_2$ the
two projections (in both cases) and by $\pi $ the fibration
$\PP (W)\to \PP ^1$, as well as the morphism $Y\to X$ given by $I
\times \pi $. Let $\Delta $ be the diagonal divisor in $X$ and
$D$ its inverse image in $Y$. Let $Q$ be the universal
quotient bundle on $\PP \Gamma (W)$ and its restriction to $\PP
(W)$. Now consider on $Y$ the bundle homomorphism
$p_1^*(W)^* \to p_2^*(Q)$ obtained as the composition of the pull
back by $p_1$ of the natural inclusion $W^* \to \Gamma
(W)^*\otimes {\cal O}$ and the pull-back by $p_2$ of the
tautological map $\Gamma (W)^* \otimes {\cal O} \to Q$. This
homomorphism is injective as a sheaf morphism but has one-dimensional
kernel on
the fibres over points of $D$. Thus we obtain an inclusion of
${\cal L}:= p_1^*(\Lambda ^{k - 1}W^*)\otimes {\cal O}(D)$ into $p_2^*(\Lambda
^{k - 1}Q)$. Note that ${\cal O}(D)$ is isomorphic to $ p_1^*{\cal O}(1)
\otimes p_2^* \pi ^*{\cal O}(1)$ so that ${\cal L}$
is isomorphic to $p_1^*({\cal O}(-k + 1))\otimes
p_2^*({\pi ^*(\cal O}(1)))$. Taking direct image by $p_1$ we get
a homomorphism of ${\cal O}(-k + 1) \otimes \Gamma (\PP ^1, {\cal
O}(1))$ into $\Gamma (\Lambda ^{k - 1}Q)$. This fits in the
following commutative diagram 

$$
\begin{array}{ccc}
{\cal O}(-k)& \rightarrow & \Lambda ^{k - 1}\Gamma (Q)\otimes {\cal O}\\
\downarrow  &             & \downarrow                  \\
{\cal O}(-k + 1)\otimes \Gamma (\PP ^1, {\cal O}(1)) & \rightarrow &
\Gamma (\Lambda ^{k - 1}Q) \otimes {\cal O} \\
\downarrow  &             & \downarrow                  \\              
{\cal O}(-k + 2)&\rightarrow & coker\Lambda ^{k - 1}\Gamma
(Q)\otimes {\cal O}\to \Gamma (\Lambda ^{k - 1}Q)\otimes {\cal O}.
\end{array}
$$

We wish to make two remarks here. Firstly the lower horizontal
arrow is nonzero. In fact, for any point $x$ of $\PP ^1$, the
middle horizontal arrow gives a two-dimensional space of
sections of $\Lambda ^{k - 1}Q$. This is obtained as follows.
Consider the sheaf inclusion of the trivial subbundle $W^*_x$ in
$\Gamma (W)^*$ on $\PP(W)$ and compose it with the natural
homomorphism of the trivial bundle with fibre $\Gamma (W)^*$
into $Q$. Take the $(k - 1)$-th exterior power of this map. This
becomes an inclusion of ${\cal O}(x) = {\cal O}(1)$ in $\Lambda ^{k -
1}Q$. Thus at the $\Gamma $-level this gives the
two-dimensional space of sections required. Clearly the sections
of the trivial bundle $\Lambda ^{k -1}(W_x^*)$ give a
one-dimensional subspace of this. This is the top horizontal
arrow in our diagram. Conversely let $s$ be an element 
of $\Lambda ^{k - 1}\Gamma (Q) = \Lambda ^{k - 1}(\Gamma W)^*$,
then its exterior product with any element of $W_x^*$ gives an
element of $\Lambda ^k\Gamma (Q) = \Lambda ^k (\Gamma W)^*$.
If $s$ is actually a section
of the sub-bundle generated by $W_x$, then this exterior product
should be zero at the generic point and hence 0.  This implies
that $s$ belongs to $\Lambda ^{k - 1}(W_x)^* = {\cal O}(-k)$.

Secondly, since all our constructions are canonical and $W$ is a
homogeneous bundle, it follows that  the lower horizontal arrow
is $SL(2)$-equivariant. Now the proposition is a consequence of
the following claim: If ${\cal O}(-n)$ admits a non-zero map
into a trivial bundle, which is equivariant for the natural
$SL(2)$-actions, then the rank of the trivial bundle is at least
$n + 1$. To prove it, use the dual map of the trivial bundle into
${\cal O}(n)$ and use the fact that $\Gamma ({\cal O}(n))$ is an
irreducible $SL(2)$-module. This implies that the induced map at
the $\Gamma $-level, which is nonzero by assumption, is actually
injective. $\hfill \square$

\section { Extra syzygies of gonal curves} \label {S5} 

We say that a curve of genus $g = 2k - 1$ is $k$-gonal if it carries
a line
bundle $L$ of degree $k$ whose linear system has no base points and thus
yields  a $k$-sheeted morphism $\pi $ onto $\PP^1$. In this paragraph,
we prove that $k$-gonal curves of genus $2k-1$
have at least $k-1$ extra syzygies.

\begin{prop} Let $C$ be a nonhyperelliptic $k$-gonal curve of
genus $2k - 1$, with $L$ the special line bundle of degree $k$ and
$Q_C$ 
the restriction of the tautological quotient bundle on $\PP^{g
-1}$ to the canonically imbedded curve $C$. Then the dimension
of $H^0(C, \Lambda ^{k - 1}Q_C)$ is at least ${g\choose k - 1} + k
- 1$.
\end{prop}

\noindent {\bf Proof.}

Consider the direct
image $V$ of the canonical line bundle $K$ of $C$ by $\pi $.  The
 so-called trace map gives a homomorphism of $V$ onto $K_{\PP
^1} = {\cal O}(-2)$. Let $W$ be its kernel. Thus we have an
exact sequence

$$
0 \rightarrow W \rightarrow V \rightarrow {\cal
O}(-2)\rightarrow 0
$$

Since ${\cal O}(-2)$ has no nonzero sections it follows that
$\Gamma (W) = \Gamma (V) = \Gamma (C,K)$.  Moreover, the kernel
of the evaluation map $\Gamma (C,K) \to W_p$ at any point $p\in
\PP ^1$ is simply the set of sections vanishing on $\pi
^{-1}(p)$, that is to say $s\Gamma (K\otimes L^{-1})$ where $s$
is a nonzero section of $L$ vanishing on this fibre. On
computing the dimension of this space to be $k$ by Riemann-Roch,
we find that the evaluation map from $\Gamma (C,K)_{\PP ^1}$ to
$V$ is actually onto $W$. Thus, $W$ is generated by global sections
and we get a morphism of $\PP (W)$ into $\PP \Gamma (C,K)$.
Finally the pull-back of $V$ to $C$, namely $\pi ^*\pi _*(K)$
comes with a natural homomorphism onto $K$. Indeed the natural
surjection of $\Gamma (K)$ onto $K$ factors through this map,
which can be thought of as `evaluation along fibres'.  Thus we
have a morphism from  $C$ to $\PP (W)$ the composition of which with
the above mentioned morphism from $\PP (W) \to \PP \Gamma(K)$
to $C$ is the canonical imbedding. Thus our claim follows from Section
\ref {scroll}.
$\hfill \square$

\section{Proof of the theorem}

In this section, we give the proof of our theorem.

We start with the

\begin{lem} Let $S$ be a smooth variety and $E$ and
$F$ two vector bundles of the same rank $n$. Let $f:E\to F$ be a
homomorphism which is generically an isomorphism, and $D$ a
subvariety of codimension 1 in $S$ on which $f$ has kernel of
rank $\geq r$, then the degeneracy divisor of $f$ contains $D$
as a component of multiplicity at least $r$.
\end{lem}

\noindent {\bf Proof.} 

Note that the question is local and
localising at the generic point of $D$, we may assume that $S$
is a discrete valuation ring with maximal ideal $\gotm M$ 
and that $f$ is a square
matrix of nonzero determinant. Then by a proper choice of basis
we may assume $f$ to be diagonal of the form $\delta_{i,j}t^{m_i}, 0
\leq
i,j \leq n$,
where $t$ is a generating parameter. Our assumption ensures
$m_i > 0$ for at least $r$ indices. Then clearly $det(f)$
is in $\gotm M^r$. 
Since the degeneracy locus is defined by $det (f)$, this
proves our assertion. $\hfill \square$

Before turning to the proof, we state again our

\begin{thm} Let $g = 2k-1 \geq 5$ be an odd integer. If the generic curve $C$
  of
  genus $g$ has no extra syzygies (i.e. 
$\Gamma (\Lambda ^{k} Q \otimes {\cal I}_C(1))=0$),
then so does any curve of genus $g$ with maximal Clifford index, namely 
$k-1.$
\end{thm}

\noindent{\bf Proof.} 

We have shown (see Section 3) that the syzygy divisor $S_g$ is the
degeneracy locus of a homomorphism of a vector bundle into
another of same rank, and (see Section 5) 
that at the generic $k$-gonal curve this
homomorphism has kernel of dimension at least $k - 1$. Thanks to the
previous
lemma, this implies that the locus of $k$-gonal curves is contained in
$S_g$ with
multiplicity at least $k - 1$.
By our computation in Section 3, the residual divisor has rational class
zero,
thus is the zero divisor (this is what we need
to know about ${\cal M}$). Thus in ${\cal M}$, curves with extra syzygies
are in the $k$-gonal divisor. Now even around curves with
automorphisms,
we can see by going to a covering where a universal curve exists, that
the locus of curves with extra syzygies is a divisor. Since the locus
of curves with automorphisms is of codimension at least two in ${\cal
  M}$, we get
that
even in ${\cal M}_g$, curves with extra syzygies
are in the $k$-gonal divisor, thus have nonmaximal
Clifford index. This proves our theorem.$\hfill \square$

\begin {rem}
We may even conclude that curves
with nonmaximal Clifford index (which have extra syzygies by [GL])
are all in the $k$-gonal divisor. Note that this result is true 
(without our assumption on the generic curve), cf  [ELMS].
\end {rem}


%
$$\br{lll}\mbox{A. Hirschowitz}&\quad\quad\quad&
\mbox{S. Ramanan}\\
\mbox{Universit\'e de Nice}&&\mbox{Tata Institute of Fundamental Research}\\
\mbox{Nice}&&\mbox{Mumbai}\\
\mbox{ah@math.unice.fr}&&\mbox{RAMANAN@math.tifr.res.in}\er$$

\end{document}